\documentclass{amsart}

\usepackage{amsmath}
\usepackage{amssymb}
\usepackage{mathrsfs}

\newtheorem{theorem}{Theorem}[section]
\newtheorem{claim}[theorem]{Claim}

\newtheorem{convention}[theorem]{Convention}
\newtheorem{observation}[theorem]{Observation}

\theoremstyle{definition}
\newtheorem{definition}[theorem]{Definition}

\newtheorem{fact}[theorem]{Fact}

\newtheorem{discussion}[theorem]{Discussion}

\newtheorem{hypothesis}[theorem]{Hypothesis}

\theoremstyle{remark}
\newtheorem{remark}[theorem]{Remark}
\newtheorem{question}[theorem]{Question}
\newtheorem{notation}[theorem]{Notation}

\newcommand{\fil}{{\rm fil}}
\newcommand{\rfil}{{\rm rfil}}

\newcommand{\otp}{{\rm otp}}

\newcommand{\Th}{{\rm Th}}

\newcommand{\tp}{{\rm tp}}

\newcommand{\lup}{{\rm l.u.p.}}
\newcommand{\lrp}{{\rm l.r.p.}}

\newcommand{\st}{{\rm st}}

\newcommand{\respects}{{\rm respects}}

\newcommand{\iif}{{\rm if}}

\newcommand{\uflp}{{\rm l.u.p.}}
\newcommand{\flp}{{\rm l.r.p.}}

\newcommand{\luft}{{\rm l.u.f.t.}}
\newcommand{\lft}{{\rm l.f.t.}}

\newcommand{\uf}{{\rm uf}}
\newcommand{\qf}{{\rm qf}}
\newcommand{\ruf}{{\rm ruf}}

\newcommand{\arity}{{\rm arity}}

\newcommand{\Mod}{{\rm Mod}}

\newcommand{\eq}{{\rm eq}}

\newcommand{\id}{{\rm id}}

\newcommand{\Dom}{{\rm Dom}}
\newcommand{\Rang}{{\rm Rang}}

\newcommand{\rest}{{\restriction}}
\newcommand{\dom}{{\rm dom}}

\newcommand{\wilog}{{\rm without loss of generality}}
\newcommand{\Wilog}{{\rm Without loss of generality}}

\newcommand{\then}{{\underline{then}}}
\newcommand{\when}{{\underline{when}}}
\newcommand{\Then}{{\underline{Then}}}

\newcommand{\Iff}{{\underline{iff}}}
\newcommand{\mn}{{\medskip\noindent}}
\newcommand{\sn}{{\smallskip\noindent}}

\newcommand{\gB}{{\mathfrak B}}

\newcommand{\cE}{{\mathscr E}}

\newcommand{\cH}{{\mathscr H}}

\newcommand{\bbL}{{\mathbb L}}

\newcommand{\bbN}{{\mathbb N}}

\newcommand{\cP}{{\mathscr P}}
\newcommand{\varp}{{\varepsilon}}

\newcommand{\cf}{{\rm cf}}

\newcount\skewfactor
\def\mathunderaccent#1#2 {\let\theaccent#1\skewfactor#2
\mathpalette\putaccentunder}
\def\putaccentunder#1#2{\oalign{$#1#2$\crcr\hidewidth
\vbox to.2ex{\hbox{$#1\skew\skewfactor\theaccent{}$}\vss}\hidewidth}}

\newenvironment{PROOF}[2][\proofname.]
   {\begin{proof}[#1]\renewcommand{\qedsymbol}{\textsquare$_{\rm #2}$}}
   {\end{proof}}

\usepackage[hidelinks]{hyperref}
\begin{document}
\makeatletter\def\shfiuwefootnote{\gdef\@thefnmark{}\@footnotetext}\makeatother\shfiuwefootnote{Version 2020-10-17. See \url{https://shelah.logic.at/papers/1101/} for possible updates.}

\title {Isomorphic limit ultrapowers for infinitary logic \\
Sh1101}
\author {Saharon Shelah}
\address{Einstein Institute of Mathematics\\
Edmond J. Safra Campus, Givat Ram\\
The Hebrew University of Jerusalem\\
Jerusalem, 9190401, Israel\\
 and \\
 Department of Mathematics\\
 Hill Center - Busch Campus \\
 Rutgers, The State University of New Jersey \\
 110 Frelinghuysen Road \\
 Piscataway, NJ 08854-8019 USA}
\email{shelah@math.huji.ac.il}
\urladdr{http://shelah.logic.at}
\thanks{The author would like to thank the Israel Science Foundation
for partial support of this research (Grant No. 1053/11).  References
like \cite[2.11=La18]{Sh:797} means we cite from \cite{Sh:797}, Claim
2.11 which has label La18, this to help if \cite{Sh:797} will be revised.
The author thanks Alice Leonhardt for the beautiful typing.  This was
separated from \cite{Sh:1019} which was first typed May 10, 2012; so
was IJM 7367}


\subjclass[2010]{Primary: 03C45; Secondary: 03C30, 03C55}

\keywords {model theory, infinitary logics, compact cardinals,
ultrapowers, ultra limits, stability, saturated models, classification
theory, isomorphic ultralimits}


\date{2020-10-17}

\begin{abstract}
The logic $\bbL^1_\theta$  
was 
introduced in \cite{Sh:797}; it is the maximal logic
below $\bbL_{\theta,\theta}$ in which a well ordering is not
definable.  We investigate it for $\theta$ a compact cardinal.  We
prove
that 
it satisfies several parallel of classical theorems on first order logic,
strengthening the thesis that it is a natural logic.  In
particular, two models are $\bbL^1_\theta$-equivalent \Iff \, for
some $\omega$-sequence of
$\theta$-complete ultrafilters, the iterated ultra-powers by it of those two
models are isomorphic.

Also for strong limit $\lambda > \theta$ of cofinality $\aleph_0$,
every complete $\bbL^1_\theta$-theory
 has a so called 
 a special model of
cardinality $\lambda$, a parallel of saturated. 
For first order theory $ T $ and singular strong limit 
cardinal $ \lambda , T $ has a so called special model 
of cardinality $ \lambda $. Using 
``special" in our context 
is justified by:
it is unique (fixing $ T $ and $ \lambda $), all
reducts of a special model are special too, so we have another proof of interpolation in
this case.
\end{abstract}

\maketitle
\numberwithin{equation}{section}
\setcounter{section}{-1}
\newpage

\centerline {Anotated Content}
\bigskip

\noindent
\S0 \quad Introduction, pg. \pageref{0}
\medskip

\noindent
\S(0A) Background and results, (label v), pg.\pageref{0A}
\medskip

\noindent
\S(0B) \quad Preliminaries, (label w,x), pg. \pageref{0B}
\bigskip

\noindent
\S1 \quad Characterizing equivalence by $\omega$-limit ultra-powers,
(label d), pg.\pageref{1}
\mn
\begin{enumerate}
\item[${{}}$]  [We characterize $\bbL^1_{< \theta}$-equivalence
  of $M_1,M_2$ by having isomorphic ultralimits by a sequence of
length $\omega$ of $\theta$-complete ultrafilters.  This logic,
  $\bbL^1_\theta$, is from \cite{Sh:797} except that here we
  restrict ourselves to $\theta$ is a compact cardinal.]
\end{enumerate}
\bigskip

\noindent
\S2 \quad Special Models, pg.\pageref{2}
\mn
\begin{enumerate}
\item[${{}}$]  [We investigate model of cardinality a strong limit
  cardinal $> \theta$ of cofinality $\aleph_0$.
We define $\lambda$-special model of complete theory $T \subseteq
  \bbL^1_\theta(\tau_T)$, for $\lambda$ as above
and prove existence and uniqueness.  We generalize some classical
theorems in model theory.]
\end{enumerate}
\newpage

\section {Introduction} \label{0}

\subsection {Background and results} \label{0A}
\bigskip

In the sixties, ultra-products were very central in model theory.
Recall Keisler \cite{Ke61}, solving the 
outstanding 
problem in model theory  of the time,
assuming an instance of GCH
characterizes elementary equivalence
in an algebraic way; that is 
by proving: 
\begin{enumerate} 
\item[$\boxplus $]
for any two models $M_1,M_2$
(of vocabulary $\tau$ of cardinality $\le \lambda$ and of cardinality\footnote{in fact ``$M_ {\ell} $ is of cardinality $ \le \lambda ^+$ suffice}
$\le \lambda$, the following are equivalent 
provided that $ 2^ \lambda = \lambda ^+$:
\item[(a)] $M_1,M_2$ are elementarily equivalent.
\item[(b)] they have isomorphic ultrapowers,
that is  
 $M^\lambda_1/D_2\approx M^\lambda_2/D_1$
 for some ultrafilter $D_{\ell} $ on a cardinal $ \lambda _ {\ell} $
\item[(c)]  $M^ \mu /D_\approx M^\mu /D$
 for some ultrafilter $D $ on some cardinal $ \mu   $
\item[(d)]  as in (c) for $ \mu = \lambda $,
\end{enumerate} 

Kochen \cite{Koc61} uses iteration on
taking ultra-powers (on a well ordered index set) to
characterize elementary equivalence.
Gaifman \cite{Gai74} uses ultra-powers on $\aleph_1$-complete
ultrafilters iterated along linear ordered index set.  Keisler
\cite{Ke63} uses general $(\aleph_0,\aleph_0)-\uflp$, see below,
Definition \ref{x16}(4) for $\kappa = \aleph_0$.   Shelah \cite{Sh:13}
proves $ \boxplus $ in ZFC, but 
with a price: we have to omit clause (d), 
and the ultrafilter is on $\mu = 2^\lambda $).

Hodges-Shelah \cite{Sh:109} is closer to the present work and see
there on earlier works, it dealt
with isomorphic ultrapowers (and isomorphic reduced powers) for the
$\theta$-complete ultrafilter (and filter) case, but note
that having isomorphic ultra-powers by
$\theta$-complete ultrafilters is not an equivalence relation.
In particular assume $\theta > \aleph_0$
is a compact cardinal and little more (we can get it by forcing over a
universe with a supercompact cardinal and a class of measurable
cardinals).  \Then \, two models have isomorphic
ultrapowers for some  $\theta$-complete ultrafilter \Iff \, in all
relevant games the isomorphism player does not lose.  Those relevant
games are of length $\zeta < \theta$ and deal with the reducts to a
sub-vocabulary of cardinality $< \theta$
and usually those games are not determined. 

The characterization \cite{Sh:109} of having isomorphic ultra-powers by
$\theta$-complete ultra-filters, necessarily
is not so ``nice" because this relation is
not an equivalence relation.  Hence having isomorphic ultra-powers
is not equivalent 
to having the same theory in some logic.

Most relevant to the present paper is \cite{Sh:797} 
which we continue here. For notational simplicity 
let $ \theta $ be an inaccessible cardinal. 
An old problem from the seventies was:

$ \boxdot $
is there a logic between $ \mathbb{L} _{\theta , {\aleph_0} }$ 
and $ \mathbb{L} _{\lambda, \theta }$ 
which satisfies interpolation?

Generally interpolation had posed  a hard problem 
in soft model theory. Another, not so precise problem 
was to find generalizations of Lindstrom theorem, see 
\cite{Van11}.  
Now \cite{Sh:797} solve the first problem 
and suggest a solution to the second problem, by
showing that the logic $ \mathbb{L} ^1_ \theta $ 
introduced there. It was proved that it satisfies 
$\boxdot $ and  give a characterization: e.g it is a maximal
logic in the interval mentioned in $ \boxdot $ 
which satisfies  non-definability of well order
in a suitable sense (see \cite[3.4=La28]{Sh:797}.

Another line of research was investigating 
infinitary logics for $ \theta $ a compact
cardinal , see \cite{Sh:1019} and history there.
We continue those two lines, investigating $ \mathbb{L} ^1_\theta $
 for $ \theta $ a compact cardinal.
 We prove that it is an interesting logic:
   it share with first order logic  several classical 
   theorems. 

We may wonder, do we have a characterization of models
being $\bbL^1_\theta$-equivalent?

In \S1 we characterize $\bbL^1_\theta$-equivalence of models by having
isomorphic iterated ultra-powers of length $\omega$.
We then 
in \S2 
prove some further generalizations
of classical model theoretic theorems, like the existence and
uniqueness of special models in $\lambda$ when $\lambda > \theta + |T|$
is strong limit of cofinality $\aleph_0$.  All this seems to
strengthen the thesis of \cite{Sh:797} that $\bbL^1_\theta$ is a
natural logic.

We thank the referee for many helpful comments.
\bigskip

Of course, success drive us to consider further
problems. For another approach 
see \cite{Sh:893}. 



\begin{question} \label{y5}
1) Assume $ \theta $ is a strong limit 
singular cardinal of cofinality $ {\aleph_0} $.  

\noindent
1) Does the logic $ \mathbb{L} _{\theta ^+ , \theta }$
restricted to models of cardinality $ \theta $ 
has interpolation?

\noindent 
2)
Is there a logic $ {\ell} $  with interpolation 
such that:
$ \mathbb{L} _{\theta ^+, \theta }
\le {\ell} \le \mathbb{l} _{\theta ^= , \theta ^+}$,
\end{question}

\begin{question}  \label{y8}
Let $ \theta $ be a compact cardinal and $ \lambda > \theta $ 
is strong limit of cofinality $ {\aleph_0} $.

\noindent
1) Does the logic $ \mathbb{L} _{\theta, \theta }$
restricted to model of cardinality $ \lambda $ 
has interpolation?

\noindent 
2) Can we characterize when a theory 
$ T \subseteq \mathbb{L}^1 _{\theta } $
of cardinality $ < \theta $ is categorical in $ \lambda  $? 

2A)
Can we then conclude that it is categorical in
other such $ \lambda $-s?

\noindent 
3) Like parts 2), 2A0 for $ T \subseteq \mathbb{L} _{\theta, \theta }$?
%
\end{question}

\subsection {Preliminaries} \label{0B}
\bigskip

\begin{hypothesis}
\label{w0}
$\theta$ is a compact uncountable cardinal (of course, we use only
restricted versions of this).
\end{hypothesis}

\begin{notation}
\label{w2}
1) Let $\varphi(\bar x)$ mean: $\varphi$ is a formula of
$\bbL_{\theta,\theta},\bar x$ is a sequence of
variables with no repetitions including the variables occurring
freely in $\varphi$ and $\ell g(\bar x) < \theta$ if not said
 otherwise.  We use $\varphi,\psi,\vartheta$ to denote formulas and
  for a statement $ \st $ let 
 $\varphi^{\st}$ or $\varphi^{[\st]}$ or $\varphi^{\iif(\st)}$
 mean 
$\varphi$ if $\st$ is true or
 1 and $\neg\varphi$ if $\st$ is false or 0.

\noindent
2) For a set $u$, usually of ordinals, let $\bar x_{[u]}
 = \langle x_\varepsilon:\varepsilon \in u \rangle$, now $u$ may be an
 ordinal but, e.g. if $u = [\alpha,\beta)$ we may write $\bar
 x_{[\alpha,\beta)}$; similarly for $\bar y_{[u]},\bar z_{[u]}$; let
 $\ell g(\bar x_{[u]}) = u$.

\noindent
3) $\tau$ denotes a vocabulary, i.e. a set of predicates and function symbols
   each with a finite number of places,
   in other words the arity $ \arity(\tau ) = {\aleph_0} $,
   see \ref{w3} on this.

\noindent
4) $T$ denotes a theory in $\bbL_{\theta,\theta}$ or $\bbL^1_\theta$
(see below); usually complete
in the vocabulary $\tau_T$ and with a model of cardinality
$\ge \theta$ if not said otherwise.

\noindent
5) Let $\Mod_T$ be the class of models of $T$.

\noindent
6) For a model $M$ let its vocabulary be $\tau_M$.
\end{notation}

\begin{remark}   \label{w3}
1) What is the problem with predicates (and function symbols)
with infinite arity? If $ \langle M_\alpha  : 
\alpha  \le \delta   \rangle , \delta $ a limit ordinal
is increasing,  even if the universe of $ M_ \delta $ 
is the union of the universes of $ M_ \alpha , \alpha < \delta $,
 this does dot determine $ M_ \delta $. 
 
 \noindent 
 2) We can still define $ \cup \{M_ \alpha : \alpha < \delta  \}  $
 by deciding $ P^{M_ \delta }= \cup \{M_ \alpha : \alpha < \delta  \} $
 for any predicate $ P $ and treating function similarly
 (so the function symbol are interpreted as partial functions
 or better deciding to use predicates only.

  Now  with care we can use 
   $ \arity (\tau ) \le \theta $ and we sometimes remark on this.
\end{remark} 

\begin{notation}
\label{w4} 
Let 
$\varepsilon,\zeta,\xi$ denote
ordinals $< \theta$.
\end{notation}

\begin{definition}
\label{w8}
1) Let $\uf_\theta(I)$ be the set of $\theta$-complete ultrafilters on
   $I$, non-principal if not said otherwise.
Let $\fil_\theta(I)$ be the set
of $\theta$-complete filters on $I$; mainly we use
$(\theta,\theta)$-regular ones (see below).

\noindent
2) $D \in \fil_\theta(I)$ is called
$(\lambda,\theta)$-regular \when \, there is a
   witness $\bar w = \langle w_t:t \in I\rangle$ which means: $w_t \in
[\lambda]^{< \theta}$ for $t \in I$ and $\alpha < \lambda \Rightarrow
   \{t:\alpha \in w_t\} \in D$.

\noindent
3) Let $\ruf_{\lambda,\theta}(I)$ be the set of $(\lambda,\theta)$-regular
$D \in \uf_\theta(I)$; let $\rfil_{\lambda,\theta}(I)$ be the
set of $(\lambda,\theta)$-regular $D \in \fil_\theta(I)$;
when $\lambda= |I|$ we may omit $\lambda$.
\end{definition}

\begin{definition}
\label{x1}
1) $\bbL_{\theta,\theta}(\tau)$ is the set of formulas of
   $\bbL_{\theta,\theta}$ in the vocabulary $\tau$.

\noindent
2) For $\tau$-models $M,N$ let $M \prec_{\bbL_{\theta,\theta}} N$
means: if $\varphi(\bar x) \in
\bbL_{\theta,\theta}(\tau_M)$ and $\bar a \in {}^{\ell g(\bar x)}M$
 then $M \models \varphi[\bar a] \Leftrightarrow N \models \varphi[\bar a]$.
\end{definition}

\noindent
And, of course
\begin{fact}
\label{x6}
For a complete $T \subseteq \bbL_{\theta,\theta}(\tau)$.

\noindent
$(\Mod_T,\prec_{\bbL_{\theta,\theta}})$ has amalgamation
and the joint embedding property (JEP), that is:
\mn
\begin{enumerate}
\item[$(a)$]  \underline{amalgamation}: if $M_0
  \prec_{\bbL_{\theta,\theta}} M_\ell$
  for $\ell=1,2$ then there are $M_3,f_1,f_2,M'_1,M'_2$ such that
\sn
\begin{enumerate}
\item[$\bullet$]  $M_0 \prec_{\bbL_{\theta,\theta}} M_3$
\sn
\item[$\bullet$]  for $\ell=1,2,f_\ell$ is a
$\prec_{\bbL_{\theta,\theta}} $-embedding of
  $M_\ell$ into $M_3$ over $M_0$, that is, for some $\tau_T$-models
  $M'_\ell$ for $\ell=1,2$ we have $M'_\ell \prec_{\bbL_{\theta,\theta}} M_3$
and $f_\ell$ is an
isomorphism from $M_\ell$ onto $M'_\ell$ over $M_0$;
\end{enumerate}
\sn
\item[$(b)$]  \underline{JEP}: if $M_1,M_2$ are
  $\bbL_{\theta,\theta}$-equivalent $\tau$-models \then \, there is a
  $\tau$-model $M_3$ and $\prec_{\bbL_{\theta,\theta}}$-embedding
$f_\ell$ of $M_\ell$ into $M_3$ for $\ell=1,2$.
\end{enumerate}
\end{fact}

\noindent
The well known generalization of \L os theorem is:
\begin{theorem}
\label{x7}
1) If $\varphi(\bar x_{[\zeta]}) \in \bbL_{\theta,\theta}(\tau),D \in
\uf_\theta(I)$ and $M_s$ is a $\tau$-model for $s \in I$ and
$f_\varepsilon \in \prod\limits_{s \in I} M_s $
for $\varepsilon <
\zeta$ \then \, $M \models
\varphi[\ldots,f_\varepsilon/D,\ldots]_{\varepsilon < \zeta}$ \Iff \,
the set $\{s \in I:M_s \models \varphi[\ldots,f_\varepsilon(s),
\ldots]_{\varepsilon < \zeta}\}$ belongs to $D$.

\noindent
2) Similarly $M \prec_{\bbL_{\theta,\theta}} M^I/D$.
\end{theorem}

\begin{definition}
\label{x17}
0) We say $X$ respects $E$ \when \, for some set $I,E$ is an equivalence relation\footnote{
Here in the interesting cases, the number of equivalent 
classes of $ E $  in infinite, and even $ \le  \theta $, pedantically not bounded by any $ \theta _* < \theta $.
}
   on $I$ and $X \subseteq I$ and $s E t \Rightarrow (s \in
 X \Leftrightarrow t \in X)$.

\noindent
1) We say $\mathbf x = (I,D,\cE)$ is a $(\kappa,\sigma)-\luft$
(limit-ultra-filter-iteration triple) \when \,:
\mn
\begin{enumerate}
\item[$(a)$]  $D$ is a filter on the set $I$
\sn
\item[$(b)$]  $\cE$ is a family of equivalence relations on $I$
\sn
\item[$(c)$]  $(\cE,\supseteq)$ is $\sigma$-directed, i.e. if
  $\alpha(*) < \sigma$ and $E_i \in \cE$ for $i < \alpha(*)$
\then \, there is $E \in \cE$
  refining $E_i$ for every $i < \alpha(*)$
\sn
\item[$(d)$]  if $E \in \cE$ then $D/E$ is a $\kappa$-complete
  ultrafilter on $I/E$ where
$D/E := \{X/E:X \in D$ and $X$ respects $E\}$.
\end{enumerate}
\mn
1A) Let $\mathbf x$ be a $(\kappa,\theta)-\lft$ mean that above we
weaken (d) to
\mn
\begin{enumerate}
\item[$(d)'$]  if $E \in \cE$ then $D/E$ is a $\kappa$-complete filter.
\end{enumerate}
\mn
2) Omitting ``$(\kappa,\sigma)$" means $(\theta,\aleph_0)$, recalling
$\theta$ is our fixed compact cardinal.

\noindent
3) Let $(I_1,D_1,\cE_1) \le^1_h (I_2,D_2,\cE_2)$ mean that:
\mn
\begin{enumerate}
\item[$(a)$]  $h$ is a function from $I_2$ onto $I_1$
\sn
\item[$(b)$]  if $E \in \cE_1$ \then \, $h^{-1} \circ E \in \cE_2$ where
  $h^{-1} \circ E = \{(s,t):s,t \in I_2$ and $h(s) Eh(t)\}$
\sn
\item[$(c)$]  if $E_1 \in \cE_1$ and $E_2 = h^{-1} \circ E_1$ \then \,
$D_1/E_1 = h''(D_2/E_2)$.
 \end{enumerate}
\end{definition}

\begin{remark}
\label{x17f}
Note that in \ref{x17}(3), if $h = \id_{I_2}$ then $I_1 = I_2$.
\end{remark}

\begin{definition}
\label{x16}
Assume $\mathbf x = (I,D,\cE)$ is a $(\kappa,\sigma)$-l.u.f.t.

\noindent
1) For a function $f$ let $\eq(f) = \{(s_1,s_2):f(s_1) = f(s_2)\}$.
   If $\bar f = \langle f_i:i < i_*\rangle$ and $i < i_* \Rightarrow
   \dom(f_i) = I$ then $\eq(\bar f) = \cap\{\eq(f_i):i < i_*\}$.

\noindent
2) For a set $U$ let $U^I|\cE = \{f \in {}^I U:\eq(f)$ is refined by some $E
\in \cE\}$.

\noindent
3) For a model $M$ let $\flp
_{\mathbf x}(M) = M^I_D|\cE = (M^I/D) \rest
   \{f/D:f \in {}^I M$ and $\eq(f)$ is refined by some $E \in \cE\}$,
pedantically (as $\arity(\tau_M)$ may be $> \aleph_0$), $M^I_D|\cE =
\cup\{M^I_D \rest E:E \in \cE\}$; $\lrp$ stands for limit reduced power.

\noindent
4) If $\mathbf x$ is $\luft$ we may in part (3) write $\uflp_{\mathbf x}(M)$.
\end{definition}

\noindent
We now give the generalization of Keisler \cite{Ke63}; Hodges-Shelah
\cite[Lemma 1,pg.80]{Sh:109} in the case $\kappa = \sigma$.
\begin{theorem}
\label{x18}
1) If $\sigma \le \kappa$ and $(I,D,\cE)$ is $(\kappa,\sigma)-\luft,\varphi =
\varphi(\bar x_{[\zeta]}) \in \bbL_{\kappa,\sigma}(\tau)$ so $\zeta < \sigma,f_\varepsilon \in
M^I|\cE$ for $\varepsilon < \zeta$ \then \, $M^I_D|\cE \models
 \varphi[\ldots,f_\varepsilon/D,\ldots]$ iff $\{s \in I:M \models
\varphi[\ldots,f_\varepsilon(s),\ldots]_{\varepsilon < \zeta}\} \in D$.

\noindent
2) Moreover $M \prec_{\bbL_{\kappa,\sigma}} M^I_D/\cE$, pedantically
   $\mathbf j = \mathbf j_{M,\mathbf x}$ is a
$\prec_{\bbL_{\kappa,\sigma}}$-elementary embedding of $M$ into
$M^I_D/\cE$ where $\mathbf j(a) = \langle a:s \in I\rangle/D$.

\noindent
3) We define $(\prod\limits_{s \in I} M_s)^I_D|\cE$ similarly when
   $\eq(\langle M_s:s \in I\rangle)$ is refined by some $E \in \cE$,
may use this more in end of the proof of \ref{d8}.
\end{theorem}

\begin{convention}
\label{x21}
 Abusing a notation 

\noindent 
1) in $\prod\limits_{s \in I} M_s/D$ we allow $f/D$
   for $f \in \prod\limits_{s \in S} M_s$ when $S \in D$.

\noindent
2) For $\bar c \in {}^\gamma(\prod\limits_{s \in I} M_s/D)$ we can find
   $\langle \bar c_s:s \in I\rangle$ such that $\bar c_s \in
{}^\gamma(M_s)$ and $\bar c =\langle \bar c_s:s \in I
\rangle/D$ which means: if $i < \ell g(\bar c)$ \then \,
   $c_{s,i} \in M_s$ and $c_i = \langle c_{s,i}:s \in I\rangle/D$.
\end{convention}

\begin{remark}
\label{x24}
1) Why the ``pedantically" in \ref{x16}(3)?  Otherwise
if $\mathbf x$ is a $(\theta,\sigma)-\luft,(\cE_{\mathbf x},
\supseteq)$ is not $\kappa^+$-directed, $\kappa < \arity(\tau)$
   then defining $\uflp_{\mathbf x}(M)$, we have freedom: if $R \in
\tau,\arity_\tau(R) \ge \kappa$, i.e. on $R^N \rest \{\bar a:\bar a
\in {}^{\arity(P)}N$ and no $E \in \cE$ refines $\eq(\bar a)\}$ so we
have no restrictions.

\noindent
2) So, e.g. for categoricity we better restrict ourselves to vocabularies
   $\tau$ such that $\arity(\tau) = \aleph_0$.
\end{remark}

\begin{definition}
\label{x28}
We say $M$ is a $\theta$-complete model \when \, for every $\varepsilon
< \theta,R_* \subseteq {}^\varepsilon M$ and $F_*:{}^\varepsilon M
\rightarrow M$ there are $R,F \in \tau_M$ such that $R^M = R_* \wedge
F^M = F_*$.
\end{definition}

\begin{observation}
\label{x31}
1) If $M$ is a $\tau$-model of cardinality $\lambda$ \then \, there is
   a $\theta$-complete expansion $M^+$ of $M$ so $\tau(M^+) \supseteq
   \tau(M)$ and $\tau(M^+)$ has cardinality $|\tau_M| +
   2^{(\|M\|^{< \theta})}$.

\noindent
2) For models $M \prec_{\bbL_{\theta,\theta}} N$ and $M^+$ as
above the following conditions are equivalent:
\mn
\begin{enumerate}
\item[$(a)$]  $N = \uflp_{\mathbf x}(M)$ identifying $a \in M$ with
  $\mathbf j_{\mathbf x}(a) \in N$, for some $(\theta,\theta)-\luft \, \mathbf x$
\sn
\item[$(b)$]  there is $N^+$ such that $M^+
  \prec_{\bbL_{\theta,\theta}} N^+$ and $N^+ \rest \tau_M$ is
  isomorphic to $N$ over $M$, in fact we can add
  $ N^+ \rest \tau _ M = N$.
\end{enumerate}
\mn
3) For a model $M$, if $(P^M,<^M)$ is a $\theta$-directed
partial order and 
$\chi = \cf(\chi) \ge \theta$ and $\lambda =
\lambda^{\|M\|} + \chi$ \then \, for some
$(\theta,\theta)-\luft \, \mathbf x$, 
the model $N := \lup_{\mathbf x}(M)$:
satisfies $(P^N,<^N)$ has a cofinal increasing sequence of
length $\chi$ and $|P^N| = \lambda$.
\end{observation}

\begin{PROOF}{\ref{x31}}
Easy, e.g.

\noindent
3) Let $M^+$ be as in part (1).
Note that $M^+$ has Skolem functions and let $T'$ be the following
set of formulas: $\Th_{\bbL_{\theta,\theta}}(M^+) \cup
\{P(x_\varp):\varp < \lambda \cdot \chi\} \cup
\{P(\sigma(x_{\varepsilon_0},\ldots,x_{\varepsilon_i},\ldots)_{i<i(*)})
\rightarrow \sigma(x_{\varepsilon_0},\dotsc,x_{\varepsilon_i},
\ldots)_{i < i(*)} < x_\varepsilon: \sigma$ is a
$\tau(M^+)$-term so $i(*) < \theta$
and $i < i(*) \Rightarrow
\varepsilon_i < \varepsilon < \lambda \cdot \chi\}$.

Clearly
\mn
\begin{enumerate}
\item[$(*)$]   $T'$ is $(< \theta)$-satisfiable in $M^+$.
\end{enumerate}
\mn
[Why?  Because if $T'' \subseteq T'$ has cardinality $< \theta$ then
the set $u = \{\varepsilon < \lambda \cdot \chi:x_\varepsilon$ appears in
$T''\}$ has cardinality $< \theta$ and let $i(*) = \otp(u)$; clearly
for each $\varepsilon \in u$ the set $\Gamma_\varepsilon = T' \cap
\{P(\sigma(x_{\varepsilon_0},\ldots)) \rightarrow
\sigma(x_{\varepsilon_0},\dotsc,x_{\varepsilon_i},\ldots)_{i<i(*)} <
x_\varepsilon:i(*) < \theta$ and $\varepsilon_i < \varepsilon$ for
$i<i(*)\}$ has cardinality $< \theta$.  Now we choose $c_\varepsilon
\in M$ by induction on $\varepsilon \in u$ such that the
assignment $x_\zeta \mapsto c_\zeta$ for $\zeta \in \varepsilon \cap
u$ in $M^+$ satisfies $\Gamma_\varepsilon$, possible because
$|\Gamma_\varepsilon| < \theta$ and $(P^M,<^M)$ is $\theta$-directed.
So the $M^+$ with the assignment $x_\varepsilon \mapsto c_\varepsilon$
for $\varepsilon \in u$ is a model of $T''$, so $T'$ is ($<
\theta)$-satisfiable indeed.]

Recalling that $|M| = \{c^{M^+}:c \in \tau(M^+)$ an individual
constant$\}$, $T'$ is realized in some
$\prec_{\bbL_{\theta,\theta}}$-elementary extension $N^+$ of $M^+$ by
the assignment $x_\varepsilon \mapsto a_\varepsilon(\varepsilon
< \lambda \cdot \chi)$.  \Wilog \, $N^+$ is the Skolem hull of
$\{a_\varepsilon:\varepsilon < \lambda \cdot \chi\}$, so $N := N^+
\rest \tau(M)$ is as required by the choice of $T'$.  Now $\mathbf x$ is
as required exists by part (2) of the claim.
\end{PROOF}

\begin{observation}
\label{x32}
1) If $\mathbf x$ is a non-trivial $(\theta,\theta)-\luft$ and $\chi =
\cf(\uflp(\theta <))$ \then \, $\chi =\chi^{< \theta}$.

\noindent
2) Also $\mu = \mu^{< \theta}$ when $\mu$ is the cardinality of
   $\lup(\theta,<)$.
\end{observation}

\begin{PROOF}{\ref{x32}}
1) By the choice of $\mathbf x$ clearly $\chi \ge \theta$.  As $\chi$ is
regular $\ge \theta$ by a theorem of Solovay \cite{So74} we have
$\chi^{< \theta} = \chi$.

\noindent
2) See the proof of \cite[2.20(3)=La27(3)]{Sh:1019}.
\end{PROOF}

We now quote \cite[Def.2.1+La8]{Sh:797}
\begin{definition}
\label{x35}
For a vocabulary $\tau,\tau$-models $M_1,M_2$, a set $\Gamma$ of
formulas in the vocabulary $\tau$ in any logic (each with
finitely many free variables if not said otherwise; see
\cite[2.9=La10(4)]{Sh:1019}),
cardinal $\theta$ and ordinal $\alpha$ we define a game $\Game =
\Game_{\Gamma,\theta,\alpha}[M_1,M_2]$
as follows, and using $(M_1,\bar b_1),(M_2,\bar b_2)$
with their natural meaning when $\Dom(\bar b_1) = \Dom(\bar b_2)$.
\mn
\begin{enumerate}
\item[$(A)$]  The moves are indexed by $n < \omega$ (but every actual play is
finite), just before the $n$-th move we have a state $\mathbf
s_n = (A^1_n,A^2_n,h^1_n,h^2_n,g_n,\beta_n,n)$
\sn
\item[$(B)$]  $\mathbf s = (A^1,A^2,h^1,h^2,g,\beta,n) = (A^1_{\mathbf
s},A^2_{\mathbf s},h^1_{\mathbf s},h^2_{\mathbf s},g_{\mathbf s},
\beta_{\mathbf s},n_{\mathbf s})$ is a state (or $n$-state or
$(\theta,n)$-state or $(\theta,< \omega)$-state) \when :
\begin{enumerate}
\item[$(a)$]   $A^\ell \in [M_\ell]^{\le \theta}$ for $\ell=1,2$
\sn
\item[$(b)$]   $\beta \le \alpha$ is an ordinal
\sn
\item[$(c)$]  $h^\ell$ is a function from $A^\ell$ into $\omega$
\sn
\item[$(d)$]   $g$ is a partial one-to-one function from $M_1$ to $M_2$ and
let $g^1_{\mathbf s} = g^1 = g_{\mathbf s} =g$ and let
$g^2_{\mathbf s} = g^2 = (g^1_{\mathbf s})^{-1}$,
\sn
\item[$(e)$]  Dom$(g^\ell) \subseteq A^\ell$ for $\ell=1,2$
\sn
\item[$(f)$]  $g$ preserves satisfaction of the formulas in $\Gamma$ and
their negations, i.e. for $\varphi(\bar x) \in \Gamma$ and $\bar a \in
{}^{\ell g(\bar x)}\text{Dom}(g)$ we have $M_1 \models \varphi[\bar
a] \Leftrightarrow M_2 \models \varphi[g(\bar a)]$
\sn
\item[$(g)$]  if $a \in \text{ Dom}(g^\ell)$ then $h^\ell(a) < n$
\end{enumerate}
\sn
\item[$(C)$]  we define the state $\mathbf s = \mathbf s_0 = \mathbf
s^0_\alpha$ by letting $n_{\mathbf s} = 0,A^1_{\mathbf s} = \emptyset =
A^2_{\mathbf s},\beta_{\mathbf s} = \alpha,h^1_{\mathbf s} = \emptyset
= h^2_{\mathbf s},g_{\mathbf s} = \emptyset$; so really $\mathbf s$ depends
only on $\alpha$
(but in general, this may not be a state for our game as possibly for
some sentence $\psi \in \Gamma$ we have
$M_1 \models \psi \Leftrightarrow M_2 \models \neg \psi$)
\sn
\item[$(D)$]  we say that a state $\mathbf t$ extends a state $\mathbf s$
when $A^\ell_{\mathbf s} \subseteq A^\ell_{\mathbf t},h^\ell_{\mathbf s}
\subseteq h^\ell_{\mathbf t}$ for $\ell=1,2$ and $g_{\mathbf s} \subseteq
g_{\mathbf t},\beta_{\mathbf s} > \beta_{\mathbf t},n_{\mathbf s} < n_{\mathbf
t}$; we say $\mathbf t$ is a
successor of $\mathbf s$ if in addition $n_{\mathbf t} = n_{\mathbf s} +1$
\sn
\item[$(E)$]  in the $n$-th move

\underline{the anti-isomorphism} player (AIS) chooses
 $(\beta_{n+1},\iota_n,A'_n)$ such that:

$\iota_n \in \{1,2\},
\beta_{n+1} < \beta_n$ and $A^{\iota_n}_n \subseteq A'_n
\in [M_{\iota_n}]^{\le \theta}$,

\underline{the isomorphism} player (ISO) chooses a state $\mathbf s_{n+1}$
such that
\begin{enumerate}
\item[$\bullet$]   $\mathbf s_{n+1}$ is a successor of $\mathbf s_n$
\sn
\item[$\bullet$]  $A^{\iota_n}_{\mathbf s_{n+1}} = A'_n$
\sn
\item[$\bullet$]  $A^{3-\iota_n}_{\mathbf s_{n+1}} =
A^{3-\iota_n}_{\mathbf s_n} \cup \Dom(g^{3-\iota_n}_{\mathbf s_{n+1}})$
\sn
\item[$\bullet$]  if $a \in A'_n \backslash A^{\iota_n}_{\mathbf s_n}$ then
$h^{\iota_n}_{\mathbf s_{n+1}}(a) \ge n+1$
\sn
\item[$\bullet$]  $\Dom(g^{\iota_n}_{\mathbf s_{n+1}}) = \{a \in
A^{\iota_n}_{\mathbf s_n}:h^{\iota_n}_{\mathbf s_n}(a) < n+1\}$ so it includes
$\Dom(g^{\iota_n}_{\mathbf s_n})$
\sn
\item[$\bullet$]  $\beta_{\mathbf s_{n+1}} = \beta_{n+1}$.
\end{enumerate}
\sn
\item[$(F)$]
\begin{itemize}
\item  the play ends when one of the players has no
legal moves (always occur
as $\beta_n < \beta_{n-1}$) and then this player loses, this
may occur for $n=0$
\item  for $\alpha = 0$ we stipulate that ISO
wins iff $\mathbf s^0_\alpha$ is a state.
\end{itemize}
\end{enumerate}
\end{definition}

\begin{definition}
\label{x38}
1) Let $\cE^{0,\tau}_{\Gamma,\theta,\alpha}$ be the class
$\{(M_1,M_2):M_1,M_2$ are $\tau$-models and in the game
$\Game_{\Gamma,\theta,\alpha}[M_1,M_2]$ the ISO player has a winning
   strategy$\}$ where $\Gamma$ is a set of formulas in the vocabulary
   $\tau$, each with finitely many free variables.

\noindent
2) $\cE^{1,\tau}_{\Gamma,\theta,\alpha}$ is the closure of
$\cE^{0,\tau}_{\Gamma,\theta,\alpha}$ to an equivalence relation (on the
class of $\tau$-models).

\noindent
3) Above we may replace $\Gamma$ by qf$(\tau)$ which means $\Gamma =$
   the set at$(\tau)$ 
   of atomic formulas
   or bs$(\tau)$ 
   of basic formulas in the vocabulary $\tau$.

\noindent
4) Above if we omit $\tau$ we mean $\tau = \tau_\Gamma$ and if we omit
   $\Gamma$ we mean bs$(\tau)$.  Abusing
   notation we may say $M_1,M_2$ are
   $\cE^{0,\tau}_{\Gamma,\theta,\alpha}$-equivalent.
\end{definition}

The following definition \ref{w41} is closely related to the beginning
of \S1, it quote \cite[Def.2.5=La13]{Sh:797} . 
\begin{definition}
\label{w41}
1) 
For a vocabulary $\tau$, the $\tau$-models $M_1,M_2$ are $\bbL^1_{<
  \theta}$-equivalent \Iff \, for every $\mu < \theta$ and $\alpha <
\mu^+$ and $\tau_1 \subseteq \tau$ of cardinality $\le \mu$, letting
$\Gamma =$ the quantifier 
free 
formulas in $\bbL(\tau)$, the models
$M_1,M_2$ are $\cE^{1,\tau_1}_{\Gamma,\mu,\alpha}$.

\noindent
2) The logic $ \mathbb{L} _{\lambda, \kappa }$ is
defined like first order logic but we allow 
conjunctions on sets  of $ < \lambda $ formulas
and we allow quantification of the form 
$ \forall \bar{ x } $ for sequences $ \bar{ x } $ 
of length $ < \kappa $; \underline{however} 
each formulas has to have $ < \kappa $ free variables;
and disjunctions and existential quantifications 
are defined naturally.

\noindent
2A) We define $ \mathbb{L} _{< \lambda , < \kappa }$
as $ \cup \{\mathbb{L} _{\lambda _1, \kappa _1}:
\lambda _1 < \lambda , \kappa _1 < \kappa \} $ ;
we may replace $ < \lambda ^+$ by $ \lambda $ 
and $ < \kappa ^+ $ by $ \kappa $.

\noindent 
3) The logic $ \mathbb{L} ^1_{\le \theta }$ is defined as follow: 
a sentence $\psi \in \bbL_{\le \theta}(\tau)$ \underline{iff} the 
sentence is defined using (or by) a triple
$(\text{qf}(\tau_1),\theta,\alpha)$ which means:
$\tau_1$ a sub-vocabulary of $\tau$ of
cardinality $\le \theta$ and $\alpha < \theta^+$ and for 
some sequence $\langle M_\beta:\beta < \beta(*)\rangle$
of $\tau_1$-models of length $\beta(*) \le \beth_{\alpha +1}(\theta)$ 
we have: $M \models \psi$ \underline{iff} $M$ is
$\cE^1_{\text{qf}(\tau_1),\theta,\alpha}$-equivalent to $M_\alpha$ for
some $\beta < \beta(*)$. 

\noindent
4) Let $\bbL^1_\kappa = \cup\{\bbL^1_{\le \theta}:\theta < \kappa\}$
   so $\bbL^1_{\theta^+} = \bbL^1_{\le \theta}$.
\end{definition}
\newpage

\section {Characterizing equivalence by $\omega$-limit ultra-powers} \label{1}

In \cite{Sh:797}, a logic $\bbL^1_{< \kappa} = \bigcup\limits_{\mu <
\kappa} \bbL^1_{\le \mu}$ is introduced
(here we consider $\kappa$ is strongly inaccessible for transparency), and is
proved to be stronger than $\bbL_{\kappa,\aleph_0}$ but weaker than
$\bbL_{\kappa,\kappa}$, has interpolation and a characterization, well
ordering not definable in it and has an addition theorem.  Also
it is the maximal logic with some such properties.

For $\kappa = \theta$, we give a characterization of when two models
are $\bbL^1_{< \theta}$-equivalent giving an additional evidence for
the logic's naturality.

\begin{convention}
\label{d2}
In this section every vocabulary $\tau$ has $\arity(\tau) =
\aleph_0$.
\end{convention}

\noindent
Recall \cite[2.11=La18]{Sh:797} which says (we expand it)
\begin{claim} 
\label{d8} 
1)
We have $M_n \equiv_{\bbL^1_{\le \theta}} M_\omega$ for $n < \omega$ 
\underline{when} clauses (b),(c) below holds and
moreover 
$M_n \models \psi[\bar a] \Leftrightarrow M_\omega \models \psi[\bar a]$
\when clauses (a)-(e) below holds, where:
\mn
\begin{enumerate}
\item[$(a)$]  $\psi(\bar z) \in \bbL^1_{\le \theta}(\tau)$ a formula
\sn
\item[$(b)$]  $M_n \prec_{\bbL_{< \partial, \theta  ^+}} M_{n+1}$ where
$\partial = \beth_{\theta  ^+}$, recalling \ref{w41}(2A)) 
\sn
\item[$(c)$]  $M_\omega := \bigcup\limits_{n < \omega} M_n$
\sn
\item[$(d)$]  $\bar a \in {}^{\ell g(\bar z)}(M_0)$
\sn
\item[$(e)$]  $\tau = \tau(M_n)$ for $n < \omega$.
\end{enumerate}

\noindent 
2) 
Assume $|\tau| \le \mu,M_n$ is a $\tau$-model and
$M_n \prec_{\bbL_{\mu^+,\mu^+}} M_{n+1}$ for $n < \omega$ and
$M_\omega = \cup\{M_n:n < \omega\}$.  \Then \, $M_0,M_\omega$ are
$\bbL^1_{\le \mu}$-equivalent.
\end{claim}

\noindent
We need two definitions before stating and proving the theorem below.
The first definition generalizes common concepts.
\begin{definition}
\label{d9}
We say that a pair of models $(M_1,M_2)$ has isomorphic
$\theta$-complete $\omega$-iterated ultrapowers \Iff \, one can find
$D_n \in \uf_\theta(I_n)$ for every $n \in \omega$ such that
$M^1_\omega \cong M^2_\omega$, when $M^\ell_\omega = \cup\{M^\ell_k:k
\in \omega\},M^\ell_0 = M_\ell$ and $M^\ell_n
\prec_{\bbL_{\theta,\theta}} (M^\ell_n)^{I_n}/D_n = M^\ell_{n+1}$ for
$\ell=1,2$ and $n <  \omega$.
\end{definition}

\noindent
For the second definition, let $\mathbf x$ be a l.u.f.t. and 
in \ref{d10} below we define ``niceness witness". 
How do we arrive to this definition?
If we try to analyze  how to prove that two 
$ \mathbb{L} ^1_ \theta $-equivalent models have 
isomorphic $ \theta $-complete $ \omega $-iterated ultrapowers
by a sequence of length $ \omega d $ of approximations;
this is a natural was to carry the induction step. 
The reader may return to this after reading the proof of 
$ (a) \rightarrow (e) $ of \ref{d11}.

To understand this (and the proof of \ref{d11} the reader may 
consider the case $ \theta = {\aleph_0} $, which naturally is simpler
and tell us that for each coordinate 
$ s \in I $ we play a game of 
an Ehrenfuecht-Fraisse game. 
Note also that Claim \ref{d8} clarify why having
$ \arity(\tau) = {\aleph_0} $ help. 

\begin{definition}
\label{d10}
If $\mathbf x = (I,D,\bar E)$ is an $\luft$ and $\bar E = \langle E_n:n \in
\omega\rangle$ \then \, $\bar w$ is a niceness witness for $(I,D,\bar
E)$ \when \,:
\mn
\begin{enumerate}
\item[(a)]  $\bar w = \langle w_{s,n},\gamma_{s,n}:s \in I,n < 
  \omega\rangle$
\sn
\item[(b)]  $w_{s,n} \subseteq \lambda_n$ and $|w_{s,n}| < \theta$ and
  $|w_{s,n}| \ge |w_{s,n+1}|$
\sn
\item[(c)]  $\gamma_{s,n} < \theta$ and $(\gamma_{s,n} >
  \gamma_{s,n+1}) \vee (\gamma_{s,n+1} = 0)$
\sn
\item[(d)]  $\gamma_{s,n} = 0 \Rightarrow w_{s,n} = \emptyset$ but
  $w_{s,0} \ne \emptyset$ and for simplicity $w_{s,0}$ is infinite for
  every $s \in I$
\sn
\item[(e)]  if $n <  \omega,u \in [\lambda_n]^{< \theta}$ then $\{s
  \in I:u \subseteq w_{s,n}\} \in D$
\sn
\item[(f)]  $w_{s,n} = w_{t,n}$ and $\gamma_{s,n} = \gamma_{t,n}$ when
  $s E_n t$.
\end{enumerate}
\end{definition}

\begin{theorem}
\label{d11}
Let $\theta$ be a compact cardinal and $M_1,M_2$ be 
two  $\tau$-models (and
$\arity(\tau) = \aleph_0$).

The following conditions are equivalent:
\mn
\begin{enumerate}
\item[$(a)$]  $M_1,M_2$ are $\bbL^1_\theta$-equivalent
\sn
\item[$(b)$]  there are $(\theta,\theta)-\luft \, \mathbf x_n = (I,D,\cE_n)$
and 
$ \cE_n \subseteq \cE_{n+1}$ for $n < \omega$ and we let $\cE =
\cup\{\cE_n:n < \omega\}$ such that
$(M_1)^I_D|\cE$ is isomorphic to $(M_2)^I_D|\cE$
\sn
\item[$(c)$]  $(M_1,M_2)$ have isomorphic $\theta$-complete
$\omega$-iterated ultrapowers, see DEf \ref{d9}
\sn
\item[$(d)$]  if $D_n \in \ruf_{\lambda_n,\theta}(I_n)$ so $|I_n| \ge
  \lambda_n$ and $\lambda_{n+1} \ge 2^{|I_n|} , 
  \lambda_n >
\|M_1\| + \|M_2\| + |\tau|$ for every $n$ \then \, the sequence
$\langle (I_n,D_n):n < \omega\rangle$ is as required in clause (c)
\sn
\item[$(e)$]  if $\mathbf x = (I,D,\cE)$ is a $\luft$ (see Def \ref{x17}(1)),
$ \cE = \{E_n:n <
  \omega\}$, for $n < \omega$ we have
$E_{n+1}$ refines $E_n,2^{|I/E_n|} \le \lambda_{n+1},
D/E_n$ is a $(\lambda_n,\theta)$-regular
$\theta$-complete ultrafilter, $\lambda_0 \ge \|M_1\| + \|M_2\| +
|\tau|,\bar w$ is a niceness witness (see Def \ref{d10}),
\then \, $\uflp_{\mathbf x}(M_1) \cong \uflp_{\mathbf x}(M_2)$, 
(see  Def. \ref{x16}(3)).
\end{enumerate}
\end{theorem}

\begin{PROOF}{\ref{d11}}
\underline{Clause $(b) \Rightarrow$ Clause $(a)$}:

So let $I,D,\cE_n(n < \omega)$ be as in clause (b) and $\cE =
\cup\{\cE_n:n < \omega\}$.  By the transitivity of being $\bbL^1_{<
  \theta}$-equivalent, clearly clause (a) follows from:
\mn
\begin{enumerate}
\item[$\boxplus_1$]  for every model $N$ the models $N,N^I_D|\cE$ are
$\bbL^1_\theta$-equivalent.
\end{enumerate}
\mn
[Why $\boxplus_1$ holds?  Let $N_n = N^I_D|\cE_n$
for $n < \omega$ and $N_\omega = \cup\{N_n:n <
\omega\}$. So by \ref{x18} we have
$N \equiv_{\bbL_{\theta,\theta}} N_0$ 
and moreover $N_n \prec_{\bbL_{\theta,\theta}} N_{n+1}$.  Hence by
\ref{d8}, that is the ``Crucial Claim"  \ref{d8} quoting
 \cite[2.11=a18]{Sh:797} we have $N_n
\equiv_{\bbL^1_{< \theta}} N_\omega$ hence 
$N \equiv_{\bbL^1_{< \theta}} N_\omega$.]
\medskip

\noindent
\underline{Clause $(c) \Rightarrow$ Clause $(b)$}:

Let $I = \prod\limits_{n < \omega} I_n,E_n = \{(\eta,\nu):\eta,\nu \in
I$ and $\eta \rest n = \nu \rest n\}$ and $D = \{X \subseteq I$: for some
$n,(\forall^{D_n} i_n \in I_n)(\forall^{D_{n-1}} i_{n-1} \in I_{n-1})
\ldots (\forall^{D_0} i_0 \in I_0)(\forall\eta)[\eta \in I \wedge
\bigwedge\limits_{\ell \le n} \eta(\ell) = i_\ell \rightarrow \eta \in
X\}$.  Now let $M^\ell_\omega \equiv (M_\ell)^I_D|\{E_n:n < \omega\}$.
 
Now 
it should be clear that
$(M_\ell)^I_D|\{E_n:n < \omega\}$ is isomorphic to $M^\ell_\omega$
 for $\ell=1,2$, so recalling $M^1_\omega \cong M^2_\omega$ by the
 present assumption, the models $(M_\ell)^I_D|\{E_n:n < \omega\}$
 for $\ell=1,2$ are isomorphic, so letting $\cE_n =
 \{E_0,\dotsc,E_n\}$ easily $(I,D,\cE_n)_{n < \omega}$ are as required
 in clause (b).
\medskip

\noindent
\underline{Clause $(d) \Rightarrow$ Clause $(c)$}:

Clause (d) is obviously stronger,
but we have to point out  that there there are such 
$ I_n, D_n, $; anyhow we shall elaborate. 
We can choose $\lambda_0 = (\|M_1\| +
\|M_2\| + |\tau| + \theta)^{< \theta},\lambda_{n+1} = 2^{\lambda_n}$ for
$n < \omega$ then letting $I_n = \lambda_n$
there is $D_n \in \ruf_{\lambda_n,\theta}(I_n)$ recalling $\theta$ is
a compact cardinal, noting $\lambda_n = \lambda_n^{< \theta}$.  Now
$\langle I_n,D_n:n < \omega\rangle$ is as required in 
the assumption of clause (d), so as we are now assuming 
clause (d), also its conclusion holds.
Now $ \langle (I_n, D_n): n < \omega \rangle $ are as required 
in   
clause (c), in
particular the isomorphism holds by 
the conclusion of 
clause (d) which as said in the previous sentence holds. .
\medskip

\noindent
\underline{Clause $(e) \Rightarrow$ Clause $(d)$}:

Let $\langle (I_n,D_n,\lambda_n):n < \omega\rangle$ be as in the
assumption of clause (d).

We define $I = \prod\limits_{n} I_n,E_n = \{(\eta,\nu):\eta,\nu \in
I,\eta \rest (n+1) = \nu \rest (n+1)\}$ and define $D$ as in the proof
of $(c) \Rightarrow (b)$ above and we choose
$\bar w = \langle w_{\eta,n}:\eta \in I,n < \omega\rangle$ as
follows.

First, choose $\bar u_n = \langle u^n_s:s \in I_n\rangle$
which witness $D_n$ is $(\lambda_n,\theta)$-regular,
i.e. $u^n_s \in [\lambda_n]^{< \theta}$
and $(\forall \alpha < \lambda_n)[\{s \in I_n:\alpha \in u^n_s\} \in D_n]$.
For $\eta \in I$ and $n < \omega$ let
$w_{\eta,n}$ be $u^n_{\eta(n)}$ if $\langle
\otp(u_{\eta(\ell)}):\ell \le n\rangle$ is decreasing and $\emptyset$
otherwise.  Let $\gamma_{\eta,n}$ be $\otp(w_{\eta,n})$.
Now we can check that the assumptions of clause (e) hold (because
of the choice of $D$), we shall elaborate two points.
First the ultrafilter $ D/ E_n $ 
is $ (\lambda, \theta  ) $-regular
because 
$ \langle \{\eta \in I: \eta ( n ) \in u^n_{\eta (n  ) } \} \rangle $
witness it.

Second 
The main point is  
to prove that $ \bar{ w } 
= \langle (w_{\eta ,n},
\gamma _{\eta, n }: \eta \in I, n < \omega \rangle $ 
is indeed a niceness witness for $ (I, D, \bar{E } )$. 
For this most clauses of \ref{d10} are easy, but we better
elaborate on clause (e) there.
that for every $n$
\mn
\begin{enumerate}
\item[$(*)_n$]  for some $X_n \in D_n$, for every $s_n \in X_n$, for
  some $X_{n-1} \in D_{n-1}$ ... for some $X_0 \in D_0$ for every $s_0
  \in X_0$, if $\langle s_0,\dotsc,s_n\rangle \trianglelefteq \eta \in
  I$ then
\sn
\begin{enumerate}
\item[(a)]  $|w_{\eta,0}| > |w_{\eta,1}| > \ldots > |w_{\eta,0}|$
\sn
\item[(b)]  $|u^\ell_{s_\ell}| > |u^{\ell +1}_{s_{\ell+1}}|$ for $\ell
  < n$.
\end{enumerate}
\end{enumerate}
\mn
Why $(*)_n$  holds?  Clause (a) holds by clause (b) and the choice of
$w_{\eta,n}$ as $u^n_{\eta(n)}$.  Clause (b) holds because
$u^{\ell+1}_{s_{\ell +1}}$ is of cardinality $< \theta$ and $\{s \in
I_\ell:|u^{\ell +1}_{s_{\ell+1}}|^+ \subseteq u^\ell_s\} \in D_\ell$.

Hence the conclusion of clause (e) holds and we are done as
in the proof of $(c) \Rightarrow (b)$.
\medskip

\noindent
\underline{Clause $(a) \Rightarrow$ Clause $(e)$}:

So assume that clause (a) holds, that is $M_1,M_2$ are
$\bbL^1_\theta$-equivalent and assume $I,D,\cE,\langle E_n:n < \omega\rangle$
and $\bar w$ are as in the assumption of clause (e), and we should
prove that its conclusion holds, that is, $\uflp_{\mathbf x}(M_1)
\cong \uflp_{\mathbf x}(M_2)$.

For every $\tau_* \subseteq \tau$ of cardinality $< \theta$ and $\mu <
\theta$, by \ref{w41} we know that $M_1 \rest \tau_*,M_2 \rest \tau_*$ are
$\bbL^1_{\le \mu}$-equivalent, hence for every $\alpha < \mu^+$ there is a
finite sequence $\langle N_{\tau_*,\mu,\alpha,k}:k \le \mathbf
k(\tau_*,\mu,\alpha)\rangle$ such that: 
\mn
\begin{enumerate}
\item[$(*)_1$]  $(a) \quad N_{\tau_*,\mu,\alpha,0} = M_1 \rest \tau_*$
\sn
\item[${{}}$]  $(b) \quad N_{\tau_*,\mu,\alpha,\mathbf
  k(\tau_*,\mu,\alpha)} = M_2 \rest \tau_*$
\sn
\item[${{}}$]  $(c) \quad$ in the game
  $\Game_{\tau_*,\mu,\alpha}[N_{\tau_*,\mu,\alpha,k},
N_{\tau_*,\mu,\alpha,k+1}]$ the ISO player has a

\hskip25pt winning strategy for each $k < \mathbf k(\tau_*,\mu,\alpha)$,
\underline{but} we stipulate

\hskip25pt a play to have $\omega$ moves, by deciding they continue to
choose the

\hskip25pt moves even when one side already wins
using the same state 
\hskip25pt except changing $n_ \mathbf{s} $
\end{enumerate} 
\sn

[Why? By Def. \ref{x35} which quote \cite[2.1=La8]{Sh:797}

\begin{enumerate} 
\item[$(*)_2$]  \wilog \, $\|N_{\tau_*,\mu,\alpha,k}\| \le \lambda_0$
for $k \in \{1,\dotsc,\mathbf k(\tau_*,\mu,\alpha)-1\}$ (even $< \theta$).
\end{enumerate}
\mn
[Why?  By (a degenerated case of ) \ref{d8}

We can (\wilog \,) assume:
\mn
\begin{enumerate}
\item[$(*)_3$]  $(a) \quad$ above
$\mathbf k(\tau_*,\mu,\alpha) =\mathbf k$
\sn 
\item[${{}}$]  $(b) \quad \tau$ has only predicates
\end{enumerate} 

\noindent 
[Why? Clause (a) 
by monotonicity in $\tau^*,\mu$ and in $\alpha$ of $M_1
\cE^{1,\tau^*}_{\qf(\tau_*),\mu,\alpha} M_2$.
Clause (b) is easy too.]

We denote:
\begin{enumerate} 
\item[$(*)_4$]  $(a) \quad \langle P_\alpha:\alpha < |\tau|\rangle$
 list the predicates of $\tau$, recall that $|\tau| \le \mu < \lambda_0$
\sn
\item[${{}}$]  $(b) \quad$ for $t \in I$ let $\tau_t = \{P_\alpha:\alpha \in
  w_{t,0} \cap |\tau|\}$
\sn
\item[$(*)_5$]  let $N_{s,k} := N_{\tau_s,|w_{s,0}|,
\gamma_{s,0}+1,k}$ for $s \in I$ and $k \le \mathbf k$.
\end{enumerate}
\mn
For $k \le \mathbf k$, let $\bar f_{k,n} = \langle f_{k,n,\alpha}:
\alpha < 2^{\lambda_n}\rangle$ list the members $f$ of
$\prod\limits_{s \in I} N_{s,k}$ such that $E_n$ refines
$\eq(f)$, so $f_{k,n,\alpha} = \langle f_{k,n,\alpha}(\eta):
\eta \in I\rangle$ but $\eta \in I \wedge \nu \in
I \wedge \eta E_n \nu \Rightarrow f_{k,n,\alpha}(\eta) =
f_{k,n,\alpha}(\nu)$.

Now
\mn
\begin{enumerate}
\item[$(*)_6$]
\begin{enumerate}
\item[(a)]   for $t \in I$ and $k < \mathbf k$ let
$\Game_{t,k}$ be the game $\Game_{\tau_t,|w_{t,0}|,\gamma_{t,0}+1}
[N_{t,k},N_{t,k+1}]$
\sn
\item[(b)]  let {\bf st}$_{t,k}$ be a winning
  strategy for the ISO player in $\Game_{t,k}$
\sn
\item[(c)]  if $t_1 E_0 t_2$ then $\langle
N_{t_\iota,k}:k \le \mathbf k\rangle$ are the
  same for $\iota=1,2$, moreover ($\Game_{t_1,k} = \Game_{t_2,k}$ and)
{\bf st}$_{t_1,k} =$ {\bf st}$_{t_2,k}$ for $k < \mathbf k$.
\end{enumerate}
\end{enumerate}
\mn
[Why clause (c)?  Because by $(*)_5$,
$N_{s,k},N_{\tau_s,|w_{s,0}|,\gamma_{s,0} +1,k}$ and $\tau_s$ depend
on $w_{s,0}$ only, so $N_{s,k}$ is determined by $(w_{s,0},k)$ hence
(by clause (e) of Th. \ref{d11} and  clause (f) from Definition \ref{d10}), $ N_{\mathbf{s}, k }$  depend just on
$(s/E_0,k)$.]

Now for each $k$ by induction on $n$ we choose
$\langle \mathbf s_{t,k,n}:t \in I\rangle$ such that:
\mn
\begin{enumerate}
\item[$(*)_7$]
\begin{enumerate}
\item[(a)]  $\mathbf s_{t,k,n}$ is a state of the game $\Game_{t,k}$
\sn
\item[(b)]  $\langle \mathbf s_{t,k,m}:m \le
  n\rangle$ is an initial segment of a play of $\Game_{t,k}$ in which
the ISO player uses the strategy {\bf st}$_{t,k}$
\sn
\item[(c)]   if $t_1 E_n t_2$ then $\mathbf s_{t_1,k,n} =
  \mathbf s_{t_2,k,n}$
\sn
\item[(d)]  $\beta_{\mathbf s_{t,k,n}} = \gamma_{t,n}$, see
Definition \ref{x35}
\sn
\item[(e)]   if $t \in I,n = \iota \mod 2$ and $\iota \in
\{0,1\}$ \then \, $A^\iota_{\mathbf s_{t,k,n}} \supseteq
\{f_{k + \iota,m,\alpha}(t):m < n$ and $\alpha \in w_{t,m}\}$, see
Definition \ref{x35}(E)
\end{enumerate}
\sn
\item[$(*)_8$]  we can carry the induction on $n$.
\end{enumerate}
\mn
[Why?  Straightforward.]
\mn
\begin{enumerate}
\item[$(*)_9$]  for each $k < \mathbf k,n < \omega,t \in I$ we define
  $h_{s,k,n}$, a partial function from $N_{s,k}$ to
  $N_{s,k+1}$ by $h_{s,k,n}(a_1) = a_2$ \Iff \, for some $m \le n,w_{s,m} \ne
\emptyset$ and $g_{\mathbf s_{t,k,m}}(a_1) = a_2$, see \ref{x35}(E).
\end{enumerate}
\mn
Now clearly:
\mn
\begin{enumerate}
\item[$\boxplus_1$]   for each $t \in I,k < \mathbf k$ and $n < \omega,
h_{s,k,n}$ is a partial one-to-one function and even a partial isomorphism
from $N_{s,k}$ to $N_{s,k+1}$, non-empty when $n > 0$ and
increasing with $n$.
\end{enumerate}
\mn
[Why?  By the choice of {\bf st}$_{t,k}$ and $(*)_7(a)$.]
\mn
\begin{enumerate}
\item[$\boxplus_2$]  let $Y_{k,n} = \{(f_1,f_2):f_\ell \in \prod\limits_{s
  \in I} \Dom(h_{s,k,n})$ for $\ell=1,2$ and $s \in I \Rightarrow
  f_2(s) = h_{s,k,n}(f_1(s))\}$
\sn
\item[$\boxplus_3$]  $\mathbf f_{k,n} = \{(f_1/D,f_2/D):(f_1,f_2) \in
  Y_{k,n}\}$ is a partial isomorphism from $M^I_1 \rest \{f/D:f \in
  \prod\limits_{s} N_{s,k}$ and $f$ respects $E_n\}$ to $M^I_2 \rest
\{f/D:f \in \prod\limits_{s} N_{s,k+1}$ and $f$ respects $E_n\}$
\sn
\item[$\boxplus_4$]  $\mathbf f_{k,n} \subseteq \mathbf
  f_{k,n+1}$
\sn
\item[$\boxplus_5$]
\begin{enumerate}
\item[(a)]   if $f_1 \in \prod\limits_{s} N_{s,k}$ and $\eq(f_1)$
is refined by $E_n$ \then \, for some $n_1 > n$ and $f_2 \in
\prod\limits_{s} N_{s,k+1}$ the pair $(f_1/D,f_2/D)$
belongs to $\mathbf f_{k,n_1}$
\sn
\item[(b)]   if $f_2 \in \prod\limits_{s} N_{s,k+1}$ and
  $\eq(f_2)$ is refined by $E_n$ then for some $n_1 > n$
and $f_1 \in \prod\limits_{s} N_{s,k}$ the pair $(f_1/D,f_2/D)$ belongs to
  $\mathbf f_{k,n_1}$.
\end{enumerate}
\end{enumerate}
\mn
[Why? By symmetry it suffices to deal with clause (a).
 For some $\alpha,f_1 = f_{k,n,\alpha}$, hence for every
$t \in \Dom(f_1),f_1(t) \in A^1_{\mathbf s_{t,k,n}}$.  We use the ``delaying
  function", $h_{\mathbf s_{t,k,n}}(f_1(t)) < \omega$ so for some $m$
the set $\{t \in I:h_{\mathbf s_{t,k,n}}(f_1(t)) \le m\}$ which
  respects $E_n$ belongs to $D$.  In particular $\{s:\gamma_{s,k,n} >
m\} \in D$, the rest should be clear recalling the regularity of each $D/E_m$.]

Letting $\cE = \{E_n:n < \omega\}$, putting together
\mn
\begin{enumerate}
\item[$(*)_{10}$]  $\mathbf f_k = \bigcup\limits_{n} \mathbf f_{k,n}$ is
  an isomorphism from $(\prod\limits_{s} N_{k,s})_D|\cE$ onto
  $(\prod\limits_{s} N_{k+1,s})_D|\cE$.
\end{enumerate}
\mn
Hence
\mn
\begin{enumerate}
\item[$(*)_{11}$]  $\mathbf f_{\mathbf k -1} \circ \ldots \circ
\mathbf f_0$ is an isomorphism from $(M_1)^I_D|\cE$ onto $(M_2)^I_D|\cE$.
\end{enumerate}
\mn
So we are done.
\end{PROOF}

\begin{discussion}
\label{d21}
1) So for our $\theta$, we get another characterization of $\bbL^1_\theta$.

\noindent
2) We may deal with universal homogeneous $(\theta,\sigma)-\uflp \,\mathbf x$,
   at least for $\sigma = \aleph_0$, using Definition \ref{x17}.
\end{discussion}

\begin{claim}
\label{d22}
In Theorem \ref{d11} if $\kappa = \kappa^{< \theta}
\ge \|M_1\| + \|M_2\|$ we can add:
\mn
\begin{enumerate}
\item[$(b)^+$]  like clause (b) of \ref{d11} but $|I| \le 2^\kappa$.
\end{enumerate}
\end{claim}

\begin{remark}
\label{d23}
Note we do not restrict $\tau = \tau(M_\ell)$.  See proof of $(*)_9$ below.
\end{remark}

\begin{PROOF}{\ref{d22}}
Clearly $(b)^+ \Rightarrow (b)$, so it is enough to prove $(b)
\Rightarrow (b)^+$; we shall assume $M_1,M_2,\kappa,\mathbf
x_n,D,\cE_n,\cE$ are as in (b) and let $g$ be an isomorphism from
$(M_1)^I_D/\cE$ onto $(M_2)^I_D/\cE$.

Let
\mn
\begin{enumerate}
\item[$(*)_1$]  $(a) \quad \cE'_n = \{E:E$ is an equivalence relation
  on $I$ with $\le \kappa$ equivalence

\hskip25pt classes such that some $E' \in \cE_n$ refines $E\}$
\sn
\item[${{}}$]  $(b) \quad$ let $\cE' = \cup\{\cE'_n:n \in \bbN\}$.
\end{enumerate}
\mn
Clearly
\mn
\begin{enumerate}
\item[$(*)_2$]  $(M_\ell)^I_D|\cE = (M_\ell)^I_D|\cE'$ for $\ell =
  1,2$.
\end{enumerate}
\mn
Let $\chi$ be large enough such that
$M_1,M_2,\kappa,D,I,\cE,\bar{\cE}' = \langle \cE'_n:n \in \bbN\rangle,g$ and
$(M_\ell)^I_D|\cE$ for $\ell=1,2$ belong to $\cH(\chi)$.  We can
choose $\gB \prec_{\bbL_{\kappa^+,\kappa^+}}(\cH(\chi),\in)$ of
cardinality $2^\kappa$ to which all the members of $\cH(\chi)$
mentioned above belong and such that $2^\kappa +1 \subseteq \gB$.  So
as $\tau = \tau(M_1) \in \gB$ and \wilog \, $|\tau| 
\le 2^{\| M_1 \|+\| M_2 \|}
\le 2^\kappa$; necessarily
$\tau \subseteq \gB$; (alternatively see the end of the  proof).
\mn
\begin{enumerate}
\item[$(*)_3$]  let
\sn
\begin{enumerate}
\item[$(a)$]  $I^* = I \cap \gB$
\sn
\item[$(b)$]  $\cE^*_n = \{E \rest I^*:E \in \cE'_n \cap \gB\}$
\sn
\item[$(c)$]  $\cE^* = \cup\{\cE^*_n:n \in \bbN\}$
\sn
\item[$(d)$]  let $D^*$ be any ultrafilter on $I^*$ which includes $\{I
  \cap I^*:I \in D \cap \gB\}$.
\end{enumerate}
\end{enumerate}
\mn
It is enough to check the following points:
\mn
\begin{enumerate}
\item[$(*)_4$]  $\mathbf x^*_n := (I^*,D^*,\cE^*_n)$ is a
 $(\theta,\theta)-\luft$ for every $n \in \omega$.
\end{enumerate}
\mn
Why?  E.g. note that if $E \in \cE^*_n$ then
for some $E' \in \cE'_n \cap \gB$ we have $E'
  \rest I^* = E$ hence $E$ has $\le \kappa$  equivalence classes.
Now for any such $E'$, as $E'$ has $\le
\kappa$-equivalence classes and belongs to $\gB$ clearly every
  $E'$-equivalence class is not disjoint to $I^*$ and every $A
  \subseteq I^*$ respecting $E$ is $A' \cap I^*$ for some $A' \in \gB$
  respecting $E'$.
So $D/E'_n,D^*/E$ are essentially equal, etc., that is, let
  $\pi_n:\cE^*_n \rightarrow \cE'_n$ be such that $E \in \cE^*_n
  \Rightarrow \pi_n(E) \rest I^* = E$ and let $\pi_{n,E}:\{A:A \subseteq
I^*$ respects $E\} \rightarrow \{A \subseteq I:A$ respects
$\pi_n(E)\}$ be such that $\pi_{n,E}(A) = B \Rightarrow B \cap I^* =
A$; in fact, those functions are uniquely determined.

So clearly $(*)_4$ follows by
\mn
\begin{enumerate}
\item[$(*)_5$]  $(a) \quad \pi_n$ is a one-to-one function from
  $\cE^*_n$ onto $\cE'_n \cap \gB$
\sn
\item[${{}}$]  $(b) \quad \pi_n$ preserves ``$E^1$ refines $E^2$" and
  its negation
\sn
\item[${{}}$]  $(c) \quad \cE^*_n$ is ($< \theta)$-directed
\sn
\item[${{}}$]  $(d) \quad$ if $n=m+1$ then $\cE^*_m \subseteq \cE^*_n$
  and $\pi_m \subseteq \pi_n$.
\end{enumerate}
\mn
Moreover
\mn
\begin{enumerate}
\item[$(*)_6$]  $(a) \quad$ if $E \in \cE^*_n$, \then \,
$\Dom(\pi_{n,E}) \subseteq \gB$ (because $2^\kappa \subseteq \gB$ is assumed)
\sn
\item[${{}}$]  $(b) \quad \pi_{n,E}$ is an isomorphism from the
Boolean Algebra $\Dom(\pi_{n,E})$ onto

\hskip25pt $\{A \subseteq I:A \,\, \respects \,\, \pi_n(E)\}$ which is
canonically isomorphic to

\hskip25pt  the Boolean Algebra
  $\cP(I/\pi_n(E))$ and also to $\cP(I^*/E)$
\sn
\item[${{}}$]  $(c) \quad D^* \cap \Dom(\pi_{n,E})$ is an ultrafilter
  which $\pi_{n,E}$ maps onto the

\hskip25pt  $D \cap \Rang(\pi_{n,E})$ which is an ultrafilter;
those ultrafilters are

\hskip25pt $\theta$-complete
\sn
\item[$(*)_7$]  $I^*$ has cardinality $\le 2^\kappa$.
\end{enumerate}
\mn
[Why?  Because $\gB$ has cardinality $\le 2^\kappa$.]
\mn
\begin{enumerate}
\item[$(*)_8$]  $(M_\ell)^{I^*}_{D^*}|\cE^*$ is isomorphic to
  $((M_\ell)^I_D|\cE') \rest \gB$ for $\ell=1,2$.
\end{enumerate}
\mn
[Why?  Let $\varkappa$ be the following function:
\mn
\begin{enumerate}
\item[$(*)_{8.1}$]  $(a) \quad \Dom(\varkappa) = (M_1)^{I_*}|\cE^*$
\sn
\item[${{}}$]  $(b) \quad$ if $f_1 \in (M_1)^{I_*}$ and $E \in
  \cE^*$ refines $\eq(f_1)$ then $f_2 := \varkappa(f_1)$ is

\hskip25pt  the unique function with domain $I$ such that $(\bigcup\limits_{n}
\pi_n)(E) \in \cE'$

\hskip25pt refines $\eq(f_2)$ and $f_2 \rest I^* = f_1$.
\end{enumerate}
\mn
Now easily $\varkappa$ induces an isomorphism as promised in $(*)_8$.]
\mn
\begin{enumerate}
\item[$(*)_9$]  $((M_1)^I_D|\cE') \rest \gB$ is isomorphic to
  $(M_2)^I_D|\cE') \rest \gB$.
\end{enumerate}
\mn
[Why?  By $(*)_2$ and the choices of $g$ (in the beginning) and of
  $\gB$ after $(*)_2$ this is obvious when $\tau = \tau(M_1)$
is included in $\gB$, which is equivalent to
  $|\tau| \le 2^\kappa$.  By recalling that
  $\arity(\tau) \le
  \aleph_0$, i.e. every predicate and function symbol of $\tau$ has
 finitely many places (see \ref{d11}), \wilog \, this holds.  That
 is, let $\tau' \subseteq \tau$ be such that for every predicate $P
 \in \tau$ there is one and only one $P' \in \tau'$ such that $\ell
 \in \{1,2\} \Rightarrow P^{M_\ell} = (P')^{M_\ell}$ and
similarly for every function symbol;
 clearly it suffices to deal with $M_1 \rest \tau',M_2 \rest \tau'$
 and $|\tau'| \le 2^{\|M_1\|} \le 2^\kappa$.]

Together we are done.
\end{PROOF}

\noindent
Note that the proof of \ref{d22} really uses $\kappa = \kappa^{<
  \theta}$, as otherwise $\cE'_n$ is not $(< \theta)$-directed.
How much is the assumption $\kappa = \kappa^{< \theta}$ needed in
\ref{d22}?
We can say something in \ref{d26}.
\begin{claim}
\label{d26}
Assume that $\kappa \ge 2^\theta$ but $\kappa^{< \theta} > \kappa$
hence for some regular $\sigma < \theta$ we have $\kappa^{< \sigma} =
\kappa < \kappa^\sigma$ and $\cf(\kappa) = \sigma$ and by \cite{So74}
we have $(\forall
\mu < \kappa)(\mu^\theta < \kappa)$; recall $\arity(\tau) = \aleph_0$.

\noindent
1) If $\langle \gB_i:i \le \sigma\rangle$ is a $\subseteq$-increasing
 continuous sequence of $\tau$-models and $\mathbf x$ is a
 $(\theta,\theta)-\luft$ \then \, $\uflp_{\mathbf x}(\gB_\sigma) =
 \cup\{\uflp_{\mathbf x}(\gB_i):i < \sigma\}$ and $i <j \Rightarrow
 \uflp_{\mathbf x}(\gB_i) \subseteq \uflp_{\mathbf x}(\gB_j)$.

\noindent
2) If $J$ is a directed partial order of cardinality $\le \sigma$ $(<
   \theta)$ and $\mathbf x_s = (I,D,\cE_s)$ is a
   $(\theta,\theta)-\luft$ for $s \in J$ such that $s <_J t
   \Rightarrow \cE_s \subseteq \cE_t$ and $M$ is a $\tau$-model \then
 \, $\uflp_{\mathbf x}(\gB) = \cup\{\uflp_{\mathbf x_s}(\gB):s \in J\}$ and
$s <_J t \Rightarrow \luft_{\mathbf x_s}(\gB) \subseteq \uflp_{\mathbf
  x_t}(\gB)$ under the natural identification.

\noindent
3) In \ref{d22}, $|I^*| \le \Sigma\{2^\partial:\partial
< \kappa\}$ is enough.
\end{claim}

\begin{PROOF}{\ref{d26}}
Straightforward.
\end{PROOF}
\newpage

\section {Special Models} \label{2}

Note that in Def. \ref{f2} below, 
$ M_n \prec _{\mathbb{L} {\theta, \theta }} M$ is not required.
The reader may in first reading ignore the 
special$^ \bullet $ case. 

\begin{definition}
\label{f2}
1)
Assume $\lambda > \theta$ is strong limit of cofinality $\aleph_0$.

We say a model $M$ is $\lambda$-special \when \, there are
$\bar\lambda,\bar M$ such that (we also may say $\bar M$ is a
$\lambda$-special sequence):
\mn
\begin{enumerate}
\item[$(a)$]  $M$ is a model of cardinality $\lambda$ with $|\tau(M)|
  < \lambda$
\sn
\item[$(b)$]
\begin{enumerate}
\item[$(\alpha)$]  $\bar \lambda = \langle \lambda_n:n \in \bbN\rangle$
\sn
\item[$(\beta)$]  $\lambda_n \le \lambda_{n+1}$
\sn
\item[$(\gamma)$]  $\theta \le \lambda_n < \lambda_{n+1} < \lambda =
  \sum\limits_{k} \lambda_k$ and stipulate $\lambda_{-1} = \theta$
\end{enumerate}
\sn
\item[$(c)$]
\begin{enumerate}
\item[$(\alpha)$]  $\bar M = \langle M_n:n < \omega\rangle$
\sn
\item[$(\beta)$]  $M_n \prec_{\bbL_{\theta,\theta}} M_{n+1}$
\sn
\item[$(\gamma)$]  $M = \bigcup\limits_{n} M_n$
\sn
\item[$(\delta)$]  $\lambda_n \ge  \|M_n\|  \ge \lambda _{n-1}$
recalling $ \lambda _{-1}= \theta $
\end{enumerate}
\sn
\item[(d)]
\begin{enumerate}
\item[$(\alpha)$]  $\bar D = \langle D_n:n \in \bbN\rangle$ and
  $\|M_n\| \le \lambda_n$
\sn
\item[$(\beta)$]  $D_n \in \ruf_{\lambda_{n-1},\theta}(\lambda_n)$
\sn
\item[$(\gamma)$]  $M^{\lambda_n}_n/D_n \prec_{\bbL_{\theta,\theta}} M_{n+1}$
  under the canonical identification (so hence $2^{\lambda_n} \le
  \lambda_{n+1}$)
\end{enumerate}
\end{enumerate}
\mn

\noindent
2) We say that the model $ M $ is $ \lambda $-special$^ \bullet $ 
when clauses (a),(b),(c) above hold but instead of clause (d) we have 
\mn
\begin{enumerate}
\item[(d)$'$]  if $\Gamma$ is an $\bbL_{\theta,\theta}$-type on $M_n$
of cardinality $\le \lambda_n$ with $\le \lambda_n$ free variables
then $\Gamma$ is realized in $M_{n+1}$.
\end{enumerate}
\end{definition}

\begin{claim}
\label{f12}
1) If for every $n < \omega$ we have
$D_n$ is a $(\lambda_n,\theta)$-regular $\theta$-complete
   ultra-filter on $I_n,M_{n+1} = (M_n)^{I_n}/D_n$ identifying $M_n$
   with its image under the canonical embedding into $M_{n+1}$ so $M_n
\prec_{\bbL_{\theta,\theta}} M_{n+1}$ and $\lambda_n \ge \|M_n\|,
\lambda = \sum\limits_{n} \lambda_n \ge \theta$ 
(equivalently $ > \theta $)
\then \,
   $\langle M_n:n \in \bbN\rangle$ is a $\lambda$-special sequence, so
   $M=  \bigcup\limits_{n} M_n$ is a $\lambda$-special model
   and $ M $  is a model of $ \Th_{\mathbb{L} ^1_ \theta }(M_1) $ .

\noindent
2) 
Assume $ \lambda > \theta , \cf( \lambda ) = {\aleph_0} $.
In Definition \ref{f2}, clause (d) indeed implies clause $(d)'$; 
so every $ \lambda $-special model/sequence 
is $ \lambda $-special$^ \bullet $ 
model/sequence.   
Also in Def. \ref{f2}, $ M $ is a model of 
$\Th_{\mathbb{L} ^1_ \theta }(M) $, in fact this follows by 
\ref{f2}(1)(d)$(\alpha ),\big( \beta ),( \gamma ) $.

\noindent 
4) Assume $ \lambda  > \theta  $ is a strong limit cardinal 
of cofinality $ {\aleph_0} $. If  $ M $ is a model of cardinality 
$ \ge \theta $ but $ < \lambda $ \then \, 
\begin{enumerate} 
\item[(A)] 
\begin{enumerate} 
\item[(a)] there is a  $ \lambda $-special sequence  $ \bar{ M }  $ 
with $ M_0 = M $ 
\item[(b)]  there is a $ \lambda $-special model $ N $ which is
a $ \prec _{\mathbb{L} ^1_\theta  }$-extension of $ M $
\item[(c)] $ \Th_{\mathbb{L} ^1_\theta (M)}$ 
has a $ \lambda $-special model.
\end{enumerate} 
\item[(B)]     
 If $ M $ is a model of cardinality $ \lambda $ \then \, 
 for some $ N, \bar{ M } , \bar{ N } $ we have:
 \begin{enumerate} 
 \item[(a)] $ \bar{ M } = \langle M_n \rangle $ is 
 a $ \lambda $-special$^ \bullet $ sequence with union $ M $
 \item[(b)]  $ \bar{ N } = \langle N_n \rangle $ is 
 a $ \lambda $-special$^ \bullet $ sequence with union $ N $ 
 \item[(c)] $ M_ n \prec _{\mathbb{L} _{\theta, \theta }} N_n$
 \end{enumerate} 
\item[(C)] If $M$ is a $\lambda$-special model and $\tau \subseteq \tau_M$,
\then \, $M \rest \tau$ is also a $\lambda$-special model.
\end{enumerate}

\noindent 
5) Assume $ \lambda > \theta > {\aleph_0} = \cf(\lambda )$.
 If $M$ is a $\lambda$-special$^\bullet $ 
 model and $\tau \subseteq \tau_M$,
\then \, $M \rest \tau$ is also a $\lambda$-special$^ \bullet $ 
model

\noindent
 6) If $ \lambda$ is strong limit $ > \theta $ of cofinality ${\aleph_0} $. 
 A model $ M $ is $ \lambda $-special iff it is $ \lambda $-special$^\bullet $.
\end{claim}

\begin{PROOF}{\ref{f12}}
1) If we assume clause (d) in Definition \ref{f2}, just by the
definition.  If we assume clause $(d)'$ in Definition \ref{f2}, use
part (2).

\noindent
2) It follows by the $(\lambda_n,\theta)$-regularity of $D_n$.

\noindent
3) Check the definition.

 \noindent 
4)   \underline{CLAUSE (A)}

We can choose an increasing sequence 
$ \langle  \lambda_n: n < \omega \rangle $ 
with limit $ \lambda $ such that $ \lambda_0 
= \| M \|^ \theta  $ 
and $ 2^ {\lambda _n }  < \lambda _{n + 1 }= 
\lambda _{n + 1} ^ \theta $. For each $ n $
we can choose a $ (\lambda , \theta  ) $-regular 
$ \theta $-complete ultrafilter on  $ \lambda _n $, 
and define $ M_n $ as in part (1).
Now use the conclusion of  part (1). 

 \underline{CLAUSE (B)}
 
 \Wilog \, the universe of $ M $ is $ \lambda $. 
 Choose $ \langle \lambda _n: n < \omega \rangle $ as
 above (except $ \lambda _0 \ge \| M \| $ of course),
 and by induction on $ n $ choose  
$ M_n  \prec _{\mathbb{L} ^1_\theta  } M $ 
of cardinality $ \lambda _n $ which include
$ \cup \{ M_k: k < n \}\cup \lambda _n $.
We now choose 
$ \langle M^*_k, M^*_  {k,n  }  : n < \omega  \rangle $
by induction on $ k $ such that 

\begin{enumerate} 
\item[(a)]  for $ k=0 $ we let $ M^*_k = M$  and 
$ M^*_{k, n }= M_n$
\item[(b)] for $ k= {\ell} + 1 $ let 
$ M^*_{k}= (M^*_{{\ell} })^{\lambda _k }/D_k$  and 
$ M^*_{k,n}= (M^*_{{\ell},n })^{\lambda _k }/D_k$.
\end{enumerate} 

There is no problem to carry the induction 
and we let $N = \cup \{M^*_{k,k } \}  $  and
$ N_ k = M^*_{k,k}$, now check

\underline{CLAUSE (C)}

Just read the definition.

\noindent 
5) Again just read the definition.

\noindent 
6) Easy too.
\end{PROOF}

\begin{remark} 
1)
In \ref{f23} we  we do not require that 
the $ \lambda _n $-s are the same and of course do not require
that the $ D_n $ are the same. Part (3) clarify this.

\noindent 
2) In Def. \ref{f2} clause (c)$( \delta ) $, 
it is enough to demand   ${}{}\lambda _n \ge 
\| M_n \| \ge \theta $.
\end{remark} 

\bigskip 

\begin{claim}
\label{f23}
1) If $\langle M^\ell_n:n \in \bbN\rangle$ is a $\lambda$-special
sequence (or just a
$ \lambda $-special$^ \bullet $ sequence)
with union $M_\ell$ for $\ell=1,2$ and
$\Th_{\bbL_{\theta,\theta}}(M^1_0) = \Th_{\bbL_{\theta,\theta}}
(M^2_0)$ \then \, $M_1,M_2$ are isomorphic.

\noindent
2) Moreover, if $n < \omega$ and $f$ is a partial function from
$M^1_n$ into $M^2_n$
which is $(M^1_n,M^2_n,\bbL_{\theta,\theta})$-elementary (i.e. $\bar a
\in {}^{\theta >}(\Dom(f)) \Rightarrow
f(\tp_{\bbL_{\theta,\theta}}(\bar a,\emptyset,M^1_n)) =
\tp_{\bbL_{\theta,\theta}}(f(\bar a),\emptyset,M^2_n))$ \then \, $f$
   can be extended to an isomorphism from $M_1$ onto $M_2$.

\noindent
3) If we weaken clause $(d)'$ of Definition \ref{f2} by weakening the
   conclusion to:  for some $k>n,\Gamma$ is realized in $M_k$ we get
an equivalent definition.


\end{claim}.

\begin{PROOF}{\ref{f23}}
1)  
By the hence and
forth argument; but we elaborate.
Let $ {\mathscr F} _n $ be the set of $ f $ such that:

\begin{enumerate} 
\item[(a)] $ f $ is a one to one function
\item[(b)] the domain of $ f $ is included in $ M^1_n $ 
\item[(c)] the range of $ f $ is included in $ M^2_n $
\item[(d)] if $ \zeta < \theta $ and $ \bar{ a } \in 
{}^{ \zeta } (M^1_n)$ and $ \bar{ b } = f( \bar{ a })  \in
{}^{ \zeta } (M^2_n) $  and $ \varphi ( \bar{ x } _{[\zeta ]}
\in \mathbb{L} _{\theta, \theta }(\tau (M_ {\ell} )$ 
\then \,  $ M^1_n \models \varphi [\bar{ a }  ] $ 
\Iff \,   
 $ M^2_n \models \varphi [\bar{ b }  ] $.
\end{enumerate} 

\medskip
 
 \noindent 
Easily 

\noindent 
\begin{enumerate} 
\item[$(*)_1$] the set $ {\mathscr F}_n $ is not empty 
\end{enumerate} 

\noindent 
[Why? Because the empty function belongs to $ {\mathscr F} _n $]

\noindent 
\begin{enumerate} 
\item[$(*)_2$] if $ f \in {\mathscr F} _n $ then 
some $ g \in {\mathscr F} _{n + 1 }$ extend $ f $ 
and $ M^1_n \subseteq \Dom(g)$
\end{enumerate} 

\noindent
[Why? By clause (d)' of  \ref{f2}(2)]

\begin{enumerate} 
\item[$(*)_3$] if $ f \in {\mathscr F} _n $ then 
some $ g \in {\mathscr F} _{n + 1 }$ extend $ f $ 
and $ M^1_n \subseteq \Rang(g)$
\end{enumerate} 

\noindent 
[Why? Similarly]

Together clearly we are done
\noindent
2) Same proof.

\noindent
3) Use suitable sub-sequences (using monotonicity).
\end{PROOF}

\noindent
Note that comparing definition \ref{f2} with the first order parallel,
in Claim \ref{f23}(1), a priory it is not given that
$\Th_{\bbL_{\theta,\theta}}(M_1) = \Th_{\bbL_{\theta,\theta}}(M_2)$
suffices.  Also \ref{f23} does not say that $\Th_{\bbL^1_\theta}(M)$
and $\lambda$ determines $M$ up to isomorphism because we demand that
$M^1_0,M^2_0$ are $\bbL^1_\theta$-equivalent.  However:
\begin{claim}
\label{f5}
Assume $\lambda > \theta$ is 
of cofinality $\aleph_0$
and    $T$ is a complete theory in $\bbL^1_\theta(\tau_T),|T| <
\lambda$ equivalently $|\tau_T| < \lambda$.

\noindent
1) If $ \lambda $ is strong limit
\Then \, $T$ has exactly one $\lambda$-special
model (up to isomorphism).

\noindent
2) 
$ T $ has at most 
one $ \lambda $-special$^ \bullet $ model of cardinality 
$ \lambda $ up to isomorphism. 
\end{claim}

\begin{PROOF}{\ref{f5}}
1) Assume $N_1,N_2$ are special models of $T$ of cardinality
 $\lambda$.  By Definition \ref{f2} for $\ell=1,2$ there is a triple
$(\bar\lambda_\ell,\bar M_\ell,\bar D_\ell)$ witnessing $N_\ell$ is
$\lambda$-special as there.

As $M_{\ell,0} \prec_{\bbL_{\theta,\theta}} M_{\ell,n}
\prec_{\bbL_{\theta,\theta}} M_{\ell,n+1} \prec_{\bbL^1_\theta}
\bigcup\limits_{m} M_{\ell,m} =
N_\ell$ for $n \in \bbN$, by \ref{x7} and \ref{d8},
we know $M_{\ell,0} \equiv_{\bbL^1_\kappa} N_\ell$, so
we can conclude that $M_{1,0}
\equiv_{\bbL^1_\kappa} M_{2,0}$ and both are models of $T$.

By \ref{d11} there is a sequence $\langle (\lambda_n,D_n):n \in
\bbN\rangle$ with $\Sigma_{n < \omega } 
\lambda_n = \lambda,2^{\lambda_n} \le
\lambda_{n+1}$ and $D_n$ a $(\lambda_n,\theta)$-regular ultrafilter on
$\lambda_n$ such that $M'_1 \cong M'_2$ when:
\mn
\begin{enumerate}
\item[$(*)$]  $M'_{\ell,0} = M_{\ell,0},M'_{\ell,n+1} =
(M'_{\ell,n})^{\lambda_n}/D_n$ and $M'_\ell = \bigcup\limits_{n}
M'_{\ell,n}$.
\end{enumerate}
\mn
So $M'_1 \cong M'_2$ by \ref{d11} and $N_1 \cong M'_1$ by \ref{f23}(1) and
$N_2 \cong M'_2$ similarly.  Together $N_1 \cong N_2$ is promised.

\noindent
2) 
The proof is similar to part of the proof of \ref{d11}
clause (a) implies clause (e).  
 i.e. by the
 hence and forth argument.
\end{PROOF}

\noindent
Now we can generalize Robinson lemma hence (see
e.g. \cite{Mak85}), gives an alternative proof
of the interpolation theorem, recall though that in \cite{Sh:797} we
do not assume the cardinal $\theta$ is compact).

\begin{claim}
\label{f8}
1) Assume $\tau_1 \cap \tau_2 = \tau_0,T_\ell$ is a complete theory in
   $\bbL^1_\theta(\tau_\ell)$ for $\ell=1,2$ and $T_0 = T_1 \cap
   T_2$.  \Then \, $T_1 \cup T_1$ has a model.

\noindent
2) We can allow in (1) the vocabularies to have more than one sort.

\noindent
3) The logic $\bbL^1_\theta$ satisfies the interpolation theorem.

\noindent
4) $\bbL^1_\theta$ has disjoint amalgamation, i.e. if $M_0
\prec_{\bbL^1_\theta} M_\ell$ for $\ell=1,2$ that is
$(M_0,c)_{c \in M_0},(M_\ell,c)_{c \in M_0}$
has the same $\bbL^1_\theta$-theory
 and $|M_1| \cap |M_2| = |M_0|$, \then \, there is $M_3$ such that
 $M_\ell \prec_{\bbL^1_\theta} M_3$ for $\ell=0,1,2$ (hence orbital
 types are well defined).

\noindent
5) $\bbL^1_\theta$ has\footnote{But the disjoint version may fail,
  e.g. if we have individual constants.} the JEP.
\end{claim}

\begin{PROOF}{\ref{f8}}
1) Let $\lambda > |\tau_1| + |\tau_2| + \theta$ be a strong limit
   cardinal of cofinality $\aleph_0$.  For $\ell = 1,2$ there is a
   $\lambda$-special model $M_\ell$ of $T_\ell$ by \ref{f12}(3).  Now
$N_\ell = M_\ell \rest \tau_0$ is a $\lambda$-special model of $T$.

By \ref{f5}(1), $N_1 \cong N_2$ so \wilog \, $N_1 = N_2$, and let $M$ be the
expansion of $N_1 = N_2$ by the predicates and functions of $M_1$ and of
$M_2$.  Clearly $M$ is a model of $T_1 \cup T_2$.

\noindent
2) Similarly.

\noindent
3) Follows as $\bbL^1_\theta$ being $\subseteq \bbL_{\theta,\theta}$
 satisfies $\theta$-compactness and part (1).

\noindent
4) Follows by (1), that is, let $\mathbf x$ be as in \ref{d11}(c) for
$M_1,M_2$.  So for every $C \subseteq M_0$ of cardinality $< \theta$,
letting $M_{C,\ell} = (M_\ell,c)_{c \in C}$ we have $N_{C,1} \cong
N_{C,2} \cong N_{C,0}$ where $N_{C,\ell} = \lup_{\mathbf
  x}(M_{C,\ell})$.  Hence $N_{C,0} \prec_{\bbL_{\theta,\theta}}
N_{C,\ell}$ for $\ell=1,2$ and we use ``$\bbL_{\theta,\theta}$ has
disjoint amalgamation".

\noindent
5) Follows by \ref{d11}.
\end{PROOF}

\begin{remark}
\label{f10}
This proof implies the generalization of preservation theorems, see
\cite{CK73}.
\end{remark}

\noindent
Recall that Ehrenfuecht-Mostowski \cite{EhMo56} aim was:  every
first order theory $T$ with infinite models has models with many automorphisms.
This fails for $\bbL_{\theta,\theta}$ and even
$\bbL_{\aleph_1,\aleph_1}$ as we can express ``$<$ is a well ordering".
What about $\bbL^1_\theta$?

\begin{claim}
\label{f22}
Assume ($\lambda,T$ are as above in \ref{f5} and) $M$ is a special model of $T$
of cardinality $\lambda$.  \Then \, $M$ has $2^\lambda$ automorphisms.
\end{claim}

\begin{PROOF}{\ref{f22}}
Let $\langle M_n:n < \omega\rangle$ witness $M$ is special.  The
result follows by the proof of \ref{f23}(2) noting that
\mn
\begin{enumerate}
\item[$(*)$]  if $f_n$ is an
  $(M_n,M_n,\bbL_{\theta,\theta}(\tau_M))$-elementary mapping \then
  \, there are $a_2 \in {}^\lambda(M_{n+1})$ and
  $f_\alpha,a_{2,\alpha} \in (M_{n+1})$ for $\alpha < \lambda_n$ such
  that
\sn
\begin{enumerate}
\item[$(a)$]  $a_{2,\alpha} \ne a_{2,\beta}$ for $\alpha < \beta <
  \lambda_n$
\sn
\item[$(b)$]  for $f_\alpha$ is an
  $(M^1_{n+1},M^2_{n+1},\bbL_{\theta,\theta}(\tau_M))$-elementary
  mapping
\sn
\item[$(c)$]  $f_\alpha \supseteq f$ and maps $a$ to $a_\alpha$.
\end{enumerate}
\end{enumerate}
\mn
Why this is possible?  Choose $a' \in M_{n+2} \backslash M_{n+1}$ and
choose $a_\alpha \in M_{n+1} \backslash \{a_\beta:\beta < \alpha\}$ by
induction on $\alpha < \lambda_n$ realizing
$\tp_{\bbL_{\theta,\theta}(\tau_T)}(a',M_n,M_{n+2})$.

Lastly, let $f_\alpha = f \cup \{(a_0,g(a_\alpha))\}$.

Why this is enough?  Should be clear.
\end{PROOF}
\newpage

\bibliographystyle{amsalpha}
\bibliography{shlhetal}
\end{document}